\documentclass[10pt]{article}
\usepackage{cite}
\usepackage{mathrsfs}
\usepackage{amsfonts}
\usepackage{amsmath}
\usepackage{amsfonts,amssymb,color}
\usepackage{dsfont}
\usepackage{curves}
\usepackage{mathrsfs}
\usepackage{pifont}
\usepackage{amssymb}
\allowdisplaybreaks

\numberwithin{equation}{section}

\date{}

\def\BigRoman{\uppercase\expandafter{\romannumeral\number\count 255 }}
\def\Romannumeral{\afterassignment\BigRoman\count255=}

\begin{document}
\title{Sufficient conditions for a graph with minimum degree to be $k$-critical with respect to $[1,b]$-odd factor
%\thanks{}
}
\author{\small  Sizhong Zhou\footnote{Corresponding
author. E-mail address: zsz\_cumt@163.com (S. Zhou)}\\
\small School of Science, Jiangsu University of Science and Technology,\\
\small Zhenjiang, Jiangsu 212100, China
}

\maketitle
\begin{abstract}
\noindent A spanning subgraph $F$ of a graph $G$ is called a $[1,b]$-odd factor if $b\equiv1$ (mod 2) and $d_F(v)\in\{1,3,\ldots,b\}$ for every
$v\in V(G)$. A graph $G$ of order $n\geq k+2$ is $k$-critical with respect to $[1,b]$-odd factor if for any $X\subseteq V(G)$ with $|X|=k$, $G-X$
has a $[1,b]$-odd factor. In this paper, we provide a size and spectral radius conditions for a graph with minimum degree to be $k$-critical with
respect to $[1,b]$-odd factor, respectively.
\\
\begin{flushleft}

(2020) Mathematics Subject Classification: 05C50, 05C70
\end{flushleft}
\end{abstract}

\section{Introduction}

Let $G=(V(G),E(G))$ be a finite, undirected and simple graph, where $V(G)$ denotes the vertex set of $G$ and $E(G)$ denotes the edge set of $G$.
The size and order of $G$ are denoted by $e(G)=|E(G)|$ and $n=|V(G)|$, respectively. Let $N_G(v)$ denote the neighbor set of a vertex $v$ in $G$.
The degree of a vertex $v$ in $G$, denoted by $d_G(v)$, is the cardinality of $N_G(v)$. The minimum degree, the number of odd components and the
number of components of $G$ are denoted by $\delta(G)$, $o(G)$ and $\omega(G)$, respectively. For a subset $S\subseteq V(G)$, let $G[S]$ denote
the subgraph of $G$ induced by $S$, and let $G-S$ denote the graph $G[V(G)-S]$. A graph $G$ of order $n$ is called $k$-connected if $n>k$ and
$G-X$ is connected for every set $X\subseteq V(G)$ with $|X|<k$, where $k$ is a positive integer. Let $K_n$ denote the complete graph of order
$n$. For two vertex-disjoint graphs $G_1$ and $G_2$, we denote by $G_1\cup G_2$ the union of $G_1$ and $G_2$, and by $G_1\vee G_2$ be the graph
formed from $G_1$ and $G_2$ by adding all the edges between $V(G_1)$ and $V(G_2)$.

Let $g$ and $f$ be two positive integer-valued functions defined on $V(G)$ such that $g(x)\leq f(x)$ holds for any $x\in V(G)$. Then a spanning
subgraph $F$ of $G$ is called a $(g,f)$-factor if $g(v)\leq d_F(v)\leq f(v)$ holds for any $v\in V(G)$. Let $a$ and $b$ be two positive integers
with $a\leq b$. A $(g,f)$-factor is called an $[a,b]$-factor if $g(v)=a$ and $f(v)=b$ for any $v\in V(G)$. A spanning subgraph $F$ of $G$ is
called a $(1,f)$-odd factor if $f(v)\equiv1$ (mod 2) and $d_F(v)\in\{1,3,\ldots,f(v)\}$ for every $v\in V(G)$. A $(1,f)$-odd factor is called a
$[1,b]$-odd factor if $f(v)=b$ is odd for any $v\in V(G)$. In fact, a perfect matching is a special $[1,b]$-odd factor when $b=1$. Let $k\geq1$
be an integer. A graph $G$ of order $n\geq k+2$ is $k$-critical with respect to $(1,f)$-odd factor if for any $X\subseteq V(G)$ with $|X|=k$,
$G-X$ contains a $(1,f)$-odd factor. A graph $G$ of order $n\geq k+2$ is $k$-critical with respect to $[1,b]$-odd factor if for any $X\subseteq V(G)$
with $|X|=k$, $G-X$ has a $[1,b]$-odd factor.

A criterion for a graph containing a perfect matching was established by Tutte \cite{Tutte} in 1947 and is one of the most important results in
factor theory. Anderson \cite{Ap} investigated the existence of perfect matchings in graphs by virtue of neighborhood condition. Egawa and Furuya
\cite{EF} obtained some results on star-free graphs having perfect matchings. Cui and Kano \cite{CK} obtained some results on the existence of
$[1,b]$-odd factors and $(1,f)$-odd factors in graphs. Topp and Vestergaard \cite{TV} discussed several conditions for a graph to contain a
$(1,f)$-odd factor. Kim, O, Park and Ree \cite{KOPR} provided an eigenvalue condition for the existence of $[1,b]$-odd factors in graphs. Favaron
\cite{Favaron} and Yu \cite{Yu} provided a characterization of $k$-factor-critical graphs, respectively. Wang and Yu \cite{WY} gave a proof of a
conjecture on $k$-factor-critical graphs. Ananchuen and Saito \cite{AS} showed some results on $k$-factor-critical graphs. Kano and Matsuda \cite{KM}
provided some results on $(1,f)$-odd factors and $k$-critical graphs with respect to $(1,f)$-odd factor. More results on graph factors were obtained
by Matsuda \cite{Mf}, Grimm, Johnsen and Shan \cite{GJS}, Haghparast and Ozeki \cite{HO}, Zhou et. al \cite{Za1,ZSL,ZPX1,ZXS,Zr}.

Let $G$ be a graph with $V(G)=\{v_1,v_2,\ldots,v_n\}$. The adjacency matrix $A(G)$ of $G$ is a matrix $(a_{ij})_{n\times n}$, where $a_{ij}=1$ if
$v_i$ is adjacent to $v_j$, and $a_{ij}=0$ otherwise. The largest eigenvalue of the adjacency matrix $A(G)$, denoted by $\rho(G)$, is called the
spectral radius of $G$.

O \cite{O} established a lower bound on the spectral radius (resp. the size) of a connected graph $G$ to guarantee that $G$ contains a perfect
matching. Zhou and Liu \cite{ZL} obtained a spectral radius condition to ensure the existence of a $[1,b]$-odd factor in a connected graph, and
they also provided a sharp upper bound for the number of edges in a connected graph without a $[1,b]$-odd factor. Zhou, Sun and Zhang \cite{ZSZ}
established a connection between spectral radius and $k$-factor-critical graphs. Zheng, Li, Luo and Wang \cite{ZLLW} provided three sufficient
conditions for a connected graph to be $k$-factor-critical. Zhou and Zhang \cite{ZZ} gave a signless Laplacian spectral radius condition for the
existence of $2k$-factor-critical graphs. More results on the relationship between spectral radius and spanning subgraphs can be found in
\cite{Wc,WZ,ZZS,ZW,Zt,ZZL,ZZL1}.

Motivated by \cite{O,ZSZ,KM} directly, we study the problem on the existence of $k$-critical graphs with respect to $[1,b]$-odd factor, and put
forward a size (resp. spectral radius) condition for the existence of $k$-critical graphs with respect to $[1,b]$-odd factor. Our main results are
shown as follows.

\medskip

\noindent{\textbf{Theorem 1.1.}} Let $b,k$ and $n$ be three positive integers with $n\equiv k$ (mod 2) and $b\equiv1$ (mod 2), and let $G$ be a
$(k+1)$-connected graph of order $n\geq\max\Big\{\frac{b^{2}k^{2}-2b^{2}k\delta-4b^{2}k-bk+b^{2}\delta^{2}+4b^{2}\delta+7b\delta+b^{2}+8b+1}{6b},
(b+5)\delta-(b+4)k-b+1+\frac{5}{b}\Big\}$ with minimum degree $\delta$. If
$$
e(G)\geq e(K_{\delta}\vee(K_{n-(b+1)\delta+bk-1}\cup(b\delta-bk+1)K_1)),
$$
then $G$ is $k$-critical with respect to $[1,b]$-odd factor, unless $G=K_{\delta}\vee(K_{n-(b+1)\delta+bk-1}\cup(b\delta-bk+1)K_1)$.

\medskip

\noindent{\textbf{Theorem 1.2.}} Let $b,k$ and $n$ be three positive integers with $n\equiv k$ (mod 2) and $b\equiv1$ (mod 2), and let $G$ be a
$(k+1)$-connected graph of order $n\geq\max\{b\delta^{2}-bk,(2b+3)\delta-bk+1\}$ with minimum degree $\delta$. If
$$
\rho(G)\geq\rho(K_{\delta}\vee(K_{n-(b+1)\delta+bk-1}\cup(b\delta-bk+1)K_1)),
$$
then $G$ is $k$-critical with respect to $[1,b]$-odd factor, unless $G=K_{\delta}\vee(K_{n-(b+1)\delta+bk-1}\cup(b\delta-bk+1)K_1)$.

\section{Some preliminaries}

Kano and Matsuda \cite{KM} provided a characterization for a graph to be $k$-critical with respect to $(1,f)$-odd factor.

\medskip

\noindent{\textbf{Lemma 2.1}} (Kano and Matsuda \cite{KM}). Let $k$ be a positive integer, and let $G$ be a graph of order $n\geq k+2$. Then $G$ is
$k$-critical with respect to $(1,f)$-odd factor if and only if
$$
o(G-S)\leq\sum\limits_{v\in S}{f(v)}-\max\Big\{\sum\limits_{v\in X}{f(v)}: X\subseteq S, |X|=k\Big\}
$$
for any $S\subseteq V(G)$ with $|S|\geq k$.

\medskip

Let $b$ be a positive odd integer. If $f(v)=b$ for any $v\in V(G)$ in Lemma 2.1, then we get the following result.

\medskip

\noindent{\textbf{Lemma 2.2.}} Let $k$ be a positive integer, and let $b$ be a positive odd integer. Then a graph $G$ of order $n\geq k+2$ is
$k$-critical with respect to $[1,b]$-odd factor if and only if
$$
o(G-S)\leq b|S|-bk
$$
for any $S\subseteq V(G)$ with $|S|\geq k$.

\medskip

\noindent{\textbf{Lemma 2.3.}} Let $k$ and $b$ be two positive integers with $b\equiv1$ (mod 2), and let $G$ be a connected graph of order $n\geq k+2$.
If $G$ is not $k$-critical with respect to $[1,b]$-odd factor, then there exists some nonempty subset $S\subseteq V(G)$ with $|S|\geq k$ such that
$$
o(G-S)>b|S|-bk.
$$
Furthermore, if $n\equiv k$ (mod 2), then
$$
o(G-S)\equiv b|S|-bk \ \ \ (\mbox{mod} \ 2)
$$
and
$$
o(G-S)\geq b|S|-bk+2.
$$

\medskip

\noindent{\it Proof.} Since $G$ is not $k$-critical with respect to $[1,b]$-odd factor, then it follows from Lemma 2.2 that
$$
o(G-S)>b|S|-bk
$$
for some nonempty subset $S\subseteq V(G)$ with $|S|\geq k$. Assume that $n\equiv k$ (mod 2). Then $o(G-S)$ is odd (resp. even) if and only if $|S|-k$
is odd (resp. even), which is equivalent to say that $b|S|-bk$ is odd (resp. even). Consequently, we obtain
$$
o(G-S)\equiv b|S|-bk \ \ \ (\mbox{mod} \ 2),
$$
and so
$$
o(G-S)\geq b|S|-bk+2.
$$
This completes the proof of Lemma 2.3. \hfill $\Box$

\medskip

\noindent{\textbf{Lemma 2.4}} (Zheng, Li, Luo and Wang \cite{ZLLW}). Let $n_1\geq n_2\geq\cdots\geq n_t\geq p\geq1$ be $(t+1)$ integers with
$n=\sum\limits_{i=1}^{t}n_i+s$ and $n_1<n-s-p(t-1)$. Then
$$
e(K_s\vee(K_{n_1}\cup K_{n_2}\cup\cdots\cup K_{n_t}))<e(K_s\vee(K_{n-s-p(t-1)}\cup(t-1)K_p)).
$$

\medskip

\noindent{\textbf{Lemma 2.5}} (Li and Feng \cite{LF}). Let $G$ be a connected graph and let $H$ be a subgraph of $G$. Then
$$
\rho(G)\geq\rho(H),
$$
where the equality holds if and only if $G=H$.

\medskip

\noindent{\textbf{Lemma 2.6}} (Fan, Lin and Lu \cite{FLL}). Let $n_1\geq n_2\geq\cdots\geq n_t\geq p\geq1$ be $(t+1)$ integers with
$n=\sum\limits_{i=1}^{t}n_i+s$ and $n_1<n-s-p(t-1)$. Then
$$
\rho(K_s\vee(K_{n_1}\cup K_{n_2}\cup\cdots\cup K_{n_t}))<\rho(K_s\vee(K_{n-s-p(t-1)}\cup(t-1)K_p)).
$$

\medskip

Let $M$ be a real $n\times n$ matrix, and let $V=\{1,2,\ldots,n\}$. Assume that $M$ can be written as
\begin{align*}
M=\left(
  \begin{array}{cccc}
    M_{11} & M_{12} & \cdots & M_{1r}\\
    M_{21} & M_{22} & \cdots & M_{2r}\\
    \vdots & \vdots & \ddots & \vdots\\
    M_{r1} & M_{r2} & \cdots & M_{rr}\\
  \end{array}
\right)
\end{align*}
with respect to the partition $\pi:V=V_1\cup V_2\cup\cdots\cup V_r$, where $M_{ij}$ denotes the submatrix (block) of $M$ obtained by rows in $V_i$ and
columns in $V_j$. Let $b_{ij}$ denote the average row sum of $M_{ij}$, that is, $b_{ij}$ is the sum of all entries in $M_{ij}$ divided by the number
of rows. Then the matrix $B(M)=(b_{ij})$ (or simply $B$) is called the quotient matrix of $M$. The partition $\pi$ is called equitable if every submatrix
$M_{ij}$ of $M$ has a constant row sum.

\medskip

\noindent{\textbf{Lemma 2.7}} (Brouwer and Haemers \cite{BH}, You, Yang, So and Xi \cite{YYSX}). Let $M$ be a real $n\times n$ matrix with an equitable
partition $\pi$, and let $B$ be the corresponding quotient matrix. Then every eigenvalue of $B$ is an eigenvalue of $M$. Furthermore, if $M$ is
nonnegative and irreducible, then the largest eigenvalue of $B$ equals the largest eigenvalue of $M$.

\medskip

The subsequent lemma is the well-known Cauchy Interlacing Theorem.

\medskip

\noindent{\textbf{Lemma 2.8}} (Haemers \cite{Hi}). Let $M$ be a Hermitian matrix of order $s$, and let $N$ be a principal submatrix of $M$ of order $t$.
Let $\lambda_1(M)\geq\lambda_2(M)\geq\cdots\geq\lambda_s(M)$ be the eigenvalues of $M$, and $\lambda_1(N)\geq\lambda_2(N)\geq\cdots\geq\lambda_t(N)$ be
the eigenvalues of $N$. Then $\lambda_i(M)\geq\lambda_i(N)\geq\lambda_{s-t+i}(M)$ for $1\leq i\leq t$.

\section{Proof of Theorem 1.1}

In this section, we verify Theorem 1.1, which provides a size condition for the existence of a $k$-critical graph with respect to $[1,b]$-odd factor.

\medskip

\noindent{\it Proof of Theorem 1.1.} Suppose to the contrary that a $(k+1)$-connected graph $G$ is not $k$-critical with respect to $[1,b]$-odd factor.
Then it follows from $n\equiv k$ (mod 2) and Lemma 2.3 that $o(G-S)\geq b|S|-bk+2$ for some nonempty subset $S\subseteq V(G)$ with $|S|\geq k$. Let
$|S|=s$ and $o(G-S)=t$. Then $t\geq bs-bk+2$ and $G$ is a spanning subgraph of $G_1=K_s\vee(K_{n_1}\cup K_{n_2}\cup\cdots\cup K_{n_{bs-bk+2}})$ for some
positive odd integers $n_1\geq n_2\geq\cdots\geq n_{bs-bk+2}$ with $\sum\limits_{i=1}^{bs-bk+2}{n_i}=n-s$. Hence, we get
\begin{align}\label{eq:3.1}
e(G)\leq e(G_1),
\end{align}
with equality if and only if $G=G_1$. Obviously, $s\geq k+1$ (otherwise, $s=k$ and so $\omega(G-S)\geq o(G-S)\geq bs-bk+2=2$, which is a contradiction
to $G$ being $(k+1)$-connected). The following proof will be divided into three cases by the value of $s$.

\noindent{\bf Case 1.} $s\geq\delta+1$.

Write $G_2=K_s\vee(K_{n-(b+1)s+bk-1}\cup(bs-bk+1)K_1)$. In view of Lemma 2.4, we conclude
\begin{align}\label{eq:3.2}
e(G_1)\leq e(G_2),
\end{align}
where the equality holds if and only if $(n_1,n_2,\ldots,n_{bs-bk+2})=(n-(b+1)s+bk-1,1,\ldots,1)$. Let
$G_*=K_{\delta}\vee(K_{n-(b+1)\delta+bk-1}\cup(b\delta-bk+1)K_1)$. Note that $n\geq(b+1)s-bk+2$. Then we obtain
\begin{align*}
e(G_2)-e(G_*)=&\binom{n-bs+bk-1}{2}+s(bs-bk+1)-\binom{n-b\delta+bk-1}{2}-\delta(b\delta-bk+1)\\
=&\frac{1}{4}(s-\delta)(-4bn+2b^{2}s+4bs+2b^{2}\delta+4b\delta-4b^{2}k-4bk+6b+4)\\
\leq&\frac{1}{4}(s-\delta)(-bn-3b((b+1)s-bk+2)+2b^{2}s+4bs+2b^{2}\delta+4b\delta-4b^{2}k-4bk+6b+4)\\
=&\frac{1}{4}(s-\delta)(-bn-b^{2}s+bs+2b^{2}\delta+4b\delta-b^{2}k-4bk+4)\\
\leq&\frac{1}{4}(s-\delta)(-bn-b^{2}(\delta+1)+b(\delta+1)+2b^{2}\delta+4b\delta-b^{2}k-4bk+4)\\
=&\frac{1}{4}(s-\delta)(-bn+b^{2}\delta+5b\delta-b^{2}k-4bk-b^{2}+b+4)\\
<&0,
\end{align*}
where the last two inequalities hold from $s\geq\delta+1$ and $n\geq(b+5)\delta-(b+4)k-b+1+\frac{5}{b}$, respectively. Thus, we deduce
$$
e(G_2)<e(G_*).
$$
Combining this with \eqref{eq:3.1} and \eqref{eq:3.2}, we conclude
$$
e(G)\leq e(G_1)\leq e(G_2)<e(G_*)=e(K_{\delta}\vee(K_{n-(b+1)\delta+bk-1}\cup(b\delta-bk+1)K_1)),
$$
which contradicts $e(G)\geq e(K_{\delta}\vee(K_{n-(b+1)\delta+bk-1}\cup(b\delta-bk+1)K_1))$.

\noindent{\bf Case 2.} $s=\delta$.

In terms of Lemma 2.4, we obtain
$$
e(G_1)\leq e(K_{\delta}\vee(K_{n-(b+1)\delta+bk-1}\cup(b\delta-bk+1)K_1)),
$$
where the equality holds if and only if $G_1=K_{\delta}\vee(K_{n-(b+1)\delta+bk-1}\cup(b\delta-bk+1)K_1)$. Together with \eqref{eq:3.1}, we deduce
$$
e(G)\leq e(K_{\delta}\vee(K_{n-(b+1)\delta+bk-1}\cup(b\delta-bk+1)K_1)),
$$
with equality holding if and only if $G=K_{\delta}\vee(K_{n-(b+1)\delta+bk-1}\cup(b\delta-bk+1)K_1)$. Observe that $K_{\delta}\vee(K_{n-(b+1)\delta+bk-1}\cup(b\delta-bk+1)K_1)$ is not $k$-critical with respect to $[1,b]$-odd factor. Consequently, we get a
contradiction.

\noindent{\bf Case 3.} $s\leq\delta-1$.

Write $G_3=K_s\vee(K_{n-s-(bs-bk+1)(\delta+1-s)}\cup(bs-bk+1)K_{\delta+1-s})$. Recall that $\delta(G)\geq\delta$ and $G$ is a spanning subgraph of
$G_1=K_s\vee(K_{n_1}\cup K_{n_2}\cup\cdots\cup K_{n_{bs-bk+2}})$ for some positive odd integers $n_1\geq n_2\geq\cdots\geq n_{bs-bk+2}$ with
$\sum\limits_{i=1}^{bs-bk+2}{n_i}=n-s$. It is easy to see that $\delta(G_1)\geq\delta$ and $n_{bs-bk+2}\geq\delta+1-s$. By virtue of Lemma 2.4, we
obtain
\begin{align}\label{eq:3.3}
e(G_1)\leq e(G_3),
\end{align}
where the equality holds if and only if $(n_1,n_2,\ldots,n_{bs-bk+2})=(n-s-(bs-bk+1)(\delta+1-s),\delta+1-s,\ldots,\delta+1-s)$. Recall that
$G_*=K_{\delta}\vee(K_{n-(b+1)\delta+bk-1}\cup(b\delta-bk+1)K_1)$. By a simple computation, we have
\begin{align}\label{eq:3.4}
e(G_3)-e(G_*)=&\binom{n-(bs-bk+1)(\delta+1-s)}{2}+(bs-bk+1)\binom{\delta+1-s}{2}\nonumber\\
&+s(bs-bk+1)(\delta+1-s)-\binom{n-b\delta+bk-1}{2}-\delta(b\delta-bk+1)\nonumber\\
=&\frac{1}{2}(s-\delta)(b^{2}s^{3}-(2b^{2}k+b^{2}\delta+2b^{2}-b)s^{2}\nonumber\\
&+(2bn+2b^{2}k\delta+b^{2}k^{2}+4b^{2}k+b^{2}-3b\delta-bk-4b)s\nonumber\\
&-(2bk+2b-2)n-b^{2}k^{2}\delta-2b^{2}k^{2}-2b^{2}k+b^{2}\delta\nonumber\\
&+3bk\delta+4bk+2b\delta-2\delta+3b-2).
\end{align}
Let $g(x)=b^{2}x^{3}-(2b^{2}k+b^{2}\delta+2b^{2}-b)x^{2}+(2bn+2b^{2}k\delta+b^{2}k^{2}+4b^{2}k+b^{2}-3b\delta-bk-4b)x-(2bk+2b-2)n
-b^{2}k^{2}\delta-2b^{2}k^{2}-2b^{2}k+b^{2}\delta+3bk\delta+4bk+2b\delta-2\delta+3b-2$ be a real function in $x$ with $x\in[k+1,\delta-1]$. By a
simple calculation, the derivative function of $g(x)$ is
$$
g'(x)=3b^{2}x^{2}-2(2b^{2}k+b^{2}\delta+2b^{2}-b)x+2bn+2b^{2}k\delta+b^{2}k^{2}+4b^{2}k+b^{2}-3b\delta-bk-4b.
$$
The symmetry axis of $g'(x)$ is $x=\frac{2bk+b\delta+2b-1}{3b}$. Since $b\geq1$ and $\delta\geq s+1\geq k+2$, we deduce
$k+1\leq\frac{2bk+b\delta+2b-1}{3b}\leq\delta-1$. When $x\in[k+1,\delta-1]$, we have
\begin{align*}
g'(x)\geq&g'\Big(\frac{2bk+b\delta+2b-1}{3b}\Big)\\
=&3b^{2}\Big(\frac{2bk+b\delta+2b-1}{3b}\Big)^{2}-2(2b^{2}k+b^{2}\delta+2b^{2}-b)\Big(\frac{2bk+b\delta+2b-1}{3b}\Big)\\
&+2bn+2b^{2}k\delta+b^{2}k^{2}+4b^{2}k+b^{2}-3b\delta-bk-4b\\
=&\frac{1}{3}(6bn-b^{2}k^{2}+2b^{2}k\delta+4b^{2}k+bk-b^{2}\delta^{2}-4b^{2}\delta-7b\delta-b^{2}-8b-1)\\
\geq&0,
\end{align*}
where the last inequality holds from $n\geq\frac{b^{2}k^{2}-2b^{2}k\delta-4b^{2}k-bk+b^{2}\delta^{2}+4b^{2}\delta+7b\delta+b^{2}+8b+1}{6b}$. It
follows from $g'(x)\geq0$ that $g(x)$ is increasing in the interval $[k+1,\delta-1]$. Combining this with $s\in[k+1,\delta-1]$ and
$n\geq(b+5)\delta-(b+4)k-b+1+\frac{5}{b}$, we get
\begin{align}\label{eq:3.5}
g(s)\geq&g(k+1)\nonumber\\
=&2n+bk-b\delta-2\delta-2\nonumber\\
\geq&2\Big((b+5)\delta-(b+4)k-b+1+\frac{5}{b}\Big)+bk-b\delta-2\delta-2\nonumber\\
=&(b+8)\delta-(b+8)k-2b+\frac{10}{b}\nonumber\\
\geq&(b+8)(k+2)-(b+8)k-2b+\frac{10}{b}\nonumber\\
=&16+\frac{10}{b}\nonumber\\
>&0.
\end{align}

In terms of \eqref{eq:3.4}, \eqref{eq:3.5} and $s\leq\delta-1$, we deduce
$$
e(G_3)<e(G_*).
$$
Together with \eqref{eq:3.1} and \eqref{eq:3.3}, we conclude
$$
e(G)\leq e(G_1)\leq e(G_3)<e(G_*)=e(K_{\delta}\vee(K_{n-(b+1)\delta+bk-1}\cup(b\delta-bk+1)K_1)),
$$
which contradicts $e(G)\geq e(K_{\delta}\vee(K_{n-(b+1)\delta+bk-1}\cup(b\delta-bk+1)K_1))$. This completes the proof of Theorem 1.1. \hfill $\Box$

\section{Proof of Theorem 1.2}

In this section, we put forward the proof of Theorem 1.2, which establishes a relationship between the adjacency spectral radius and $k$-critical graph
with respect to $[1,b]$-odd factor.

\medskip

\noindent{\it Proof of Theorem 1.2.} Suppose to the contrary that a $(k+1)$-connected graph $G$ is not $k$-critical with respect to $[1,b]$-odd factor.
Then by $n\equiv k$ (mod 2) and Lemma 2.3, we have $o(G-S)\geq b|S|-bk+2$ for some nonempty subset $S\subseteq V(G)$ with $|S|\geq k$. Let $|S|=s$ and
$o(G-S)=t$. Then $t\geq bs-bk+2$ and $G$ is a spanning subgraph of $G_1=K_s\vee(K_{n_1}\cup K_{n_2}\cup\cdots\cup K_{n_{bs-bk+2}})$ for some positive
odd integers $n_1\geq n_2\geq\cdots\geq n_{bs-bk+2}$ with $\sum\limits_{i=1}^{bs-bk+2}{n_i}=n-s$. In terms of Lemma 2.5, we conclude
\begin{align}\label{eq:4.1}
\rho(G)\leq\rho(G_1),
\end{align}
with equality if and only if $G=G_1$. Since $G$ is $(k+1)$-connected, we deduce $\delta\geq k+1$.

\medskip

\noindent{\bf Claim 1.} $s\geq k+1$.

\noindent{\it Proof.} Assume that $s=k$. Then $\omega(G-S)\geq o(G-S)\geq bs-bk+2=2$, which is a contradiction to $G$ being $(k+1)$-connected. This
completes the proof of Claim 1. \hfill $\Box$

\medskip

The following proof will be divided into three possible cases by the value of $s$.

\noindent{\bf Case 1.} $s\geq\delta+1$.

Let $G_2=K_s\vee(K_{n-(b+1)s+bk-1}\cup(bs-bk+1)K_1)$. By virtue of Lemma 2.6, we deduce
\begin{align}\label{eq:4.2}
\rho(G_1)\leq\rho(G_2),
\end{align}
with equality holding if and only if $(n_1,n_2,\ldots,n_{bs-bk+2})=(n-(b+1)s+bk-1,1,\ldots,1)$.

The quotient matrix of $A(G_2)$ corresponding to the partition $V(G_2)=V(K_s)\cup V(K_{n-(b+1)s+bk-1})\cup V((bs-bk+1)K_1)$ is given as
\begin{align*}
B_2=\left(
  \begin{array}{ccc}
    s-1 & n-(b+1)s+bk-1 & bs-bk+1\\
    s & n-(b+1)s+bk-2 & 0\\
    s & 0 & 0\\
  \end{array}
\right).
\end{align*}
Then the characteristic polynomial of $B_2$ is
\begin{align*}
\varphi_{B_2}(x)=&x^{3}+(-n+bs-bk+3)x^{2}+(-n-bs^{2}+bks+bs-s-bk+2)x\\
&-b(b+1)s^{3}+bns^{2}+2b^{2}ks^{2}+bks^{2}-(3b+1)s^{2}-bkns+ns-b^{2}k^{2}s+3bks-2s.
\end{align*}
Since the partition $V(G_2)=V(K_s)\cup V(K_{n-(b+1)s+bk-1})\cup V((bs-bk+1)K_1)$ is equitable, it follows from Lemma 2.7 that the largest root,
say $\theta_1$, of $\varphi_{B_2}(x)=0$ is equal to $\rho(G_2)$, that is, $\rho(G_2)=\theta_1$. Let $\theta_1=\rho(G_2)\geq\theta_2\geq\theta_3$
be the three roots of $\varphi_{B_2}(x)=0$ and $Q_2=\mbox{diag}(s,n-(b+1)s+bk-1,bs-bk+1)$. One easily checks that
\begin{align*}
Q_2^{\frac{1}{2}}B_2Q_2^{-\frac{1}{2}}=\left(
  \begin{array}{ccc}
    s-1 & s^{\frac{1}{2}}(n-(b+1)s+bk-1)^{\frac{1}{2}} & s^{\frac{1}{2}}(bs-bk+1)^{\frac{1}{2}}\\
    s^{\frac{1}{2}}(n-(b+1)s+bk-1)^{\frac{1}{2}} & n-(b+1)s+bk-2 & 0\\
    s^{\frac{1}{2}}(bs-bk+1)^{\frac{1}{2}} & 0 & 0\\
  \end{array}
\right)
\end{align*}
is symmetric, and also contains
\begin{align*}
\left(
  \begin{array}{ccc}
    n-(b+1)s+bk-2 & 0\\
    0 & 0\\
  \end{array}
\right)
\end{align*}
as its submatrix. Since $Q_2^{\frac{1}{2}}B_2Q_2^{-\frac{1}{2}}$ and $B_2$ admit the same eigenvalues, it follows from Lemma 2.8 and $s\geq\delta+1$
that
\begin{align}\label{eq:4.3}
\theta_2\leq n-(b+1)s+bk-2\leq n-(b+1)(\delta+1)+bk-2.
\end{align}

Let $G_*=K_{\delta}\vee(K_{n-(b+1)\delta+bk-1}\cup(b\delta-bk+1)K_1)$. Then the quotient matrix of $A(G_*)$ corresponding to the partition
$V(G_*)=V(K_{\delta})\cup V(K_{n-(b+1)\delta+bk-1})\cup V((b\delta-bk+1)K_1)$ equals
\begin{align*}
B_*=\left(
  \begin{array}{ccc}
    \delta-1 & n-(b+1)\delta+bk-1 & b\delta-bk+1\\
    \delta & n-(b+1)\delta+bk-2 & 0\\
    \delta & 0 & 0\\
  \end{array}
\right),
\end{align*}
whose characteristic polynomial is
\begin{align*}
\varphi_{B_*}(x)=&x^{3}+(-n+b\delta-bk+3)x^{2}+(-n-b\delta^{2}+bk\delta+b\delta-\delta-bk+2)x\\
&-b(b+1)\delta^{3}+bn\delta^{2}+2b^{2}k\delta^{2}+bk\delta^{2}-(3b+1)\delta^{2}-bkn\delta+n\delta-b^{2}k^{2}\delta+3bk\delta-2\delta.
\end{align*}
One easily see that the partition $V(G_*)=V(K_{\delta})\cup V(K_{n-(b+1)\delta+bk-1})\cup V((b\delta-bk+1)K_1)$ is equitable. By means of Lemma 2.7,
the largest root, say $\theta_*$, of $\varphi_{B_*}(x)=0$ is equal to $\rho(G_*)$, that is, $\rho(G_*)=\theta_*$. Since $K_{n-b\delta+bk-1}$ is a
proper subgraph of $G_*$, it follows from \eqref{eq:4.3} and Lemma 2.5 that
\begin{align}\label{eq:4.4}
\theta_*=\rho(G_*)>\rho(K_{n-b\delta+bk-1})=n-b\delta+bk-2>n-(b+1)(\delta+1)+bk-2\geq\theta_2.
\end{align}
Note that $\varphi_{B_*}(\theta_*)=0$. We are to prove $\varphi_{B_2}(\theta_*)=\varphi_{B_2}(\theta_*)-\varphi_{B_*}(\theta_*)>0$. By a direct
calculation, we have
\begin{align}\label{eq:4.5}
\varphi_{B_2}(\theta_*)=\varphi_{B_2}(\theta_*)-\varphi_{B_*}(\theta_*)=(s-\delta)f_1(\theta_*),
\end{align}
where $f_1(\theta_*)=b\theta_*^{2}+(-bs-b\delta+bk+b-1)\theta_*-b(b+1)(s^{2}+\delta s+\delta^{2})+(s+\delta)(bn+2b^{2}k+bk-3b-1)-bkn+n-b^{2}k^{2}+3bk-2$.
The symmetry axis of $f_1(\theta_*)$ is $\theta_*=\frac{bs+b\delta-bk-b+1}{2b}$, which implies that $f_1(\theta_*)$ is increasing in the interval
$[\frac{bs+b\delta-bk-b+1}{2b},+\infty)$. According to \eqref{eq:4.4}, $s\geq\delta+1$ and $n\geq(b+1)s-bk+2$, we easily deduce that
$$
\frac{bs+b\delta-bk-b+1}{2b}<n-b\delta+bk-2<\theta_*,
$$
and so
\begin{align}\label{eq:4.6}
f_1(\theta_*)>&f_1(n-b\delta+bk-2)\nonumber\\
=&b(n-b\delta+bk-2)^{2}+(-bs-b\delta+bk+b-1)(n-b\delta+bk-2)\nonumber\\
&-b(b+1)(s^{2}+\delta s+\delta^{2})+(s+\delta)(bn+2b^{2}k+bk-3b-1)\nonumber\\
&-bkn+n-b^{2}k^{2}+3bk-2\nonumber\\
=&-b(b+1)s^{2}+(b^{2}k-b\delta+bk-b-1)s+bn^{2}+(-2b^{2}\delta+2b^{2}k-3b)n\nonumber\\
&+(b^{3}-b)\delta^{2}+(-2b^{3}k+3b^{2}+bk-1)\delta+b^{3}k^{2}-3b^{2}k+2b.
\end{align}
Let $f_2(s)=-b(b+1)s^{2}+(b^{2}k-b\delta+bk-b-1)s+bn^{2}+(-2b^{2}\delta+2b^{2}k-3b)n+(b^{3}-b)\delta^{2}+(-2b^{3}k+3b^{2}+bk-1)\delta+b^{3}k^{2}-3b^{2}k+2b$.
Notice that
$$
\frac{b^{2}k-b\delta+bk-b-1}{2b(b+1)}<\delta+1\leq s\leq\frac{n+bk-2}{b+1}
$$
by $\delta\geq k+1$, $s\geq\delta+1$ and $n\geq(b+1)s-bk+2$. Thus, we deduce
\begin{align*}
f_2(s)\geq&f_2\Big(\frac{n+bk-2}{b+1}\Big)\\
=&-b(b+1)\Big(\frac{n+bk-2}{b+1}\Big)^{2}+(b^{2}k-b\delta+bk-b-1)\Big(\frac{n+bk-2}{b+1}\Big)\\
&+bn^{2}+(-2b^{2}\delta+2b^{2}k-3b)n+(b^{3}-b)\delta^{2}+(-2b^{3}k+3b^{2}+bk-1)\delta\\
&+b^{3}k^{2}-3b^{2}k+2b\\
=&\frac{b^{2}n^{2}}{b+1}+\frac{(-(2b^{3}+2b^{2}+b)\delta+2b^{3}k+b^{2}k+bk-3b^{2}-1)n}{b+1}\\
&+(b^{3}-b)\delta^{2}+\frac{(-2b^{4}k-2b^{3}k+bk+3b^{3}+3b^{2}+b-1)\delta}{b+1}\\
&+\frac{b^{4}k^{2}+b^{3}k^{2}+b^{2}k^{2}-3b^{3}k-2b^{2}k-3bk+2b^{2}+2}{b+1}\\
\geq&\frac{b^{2}((2b+3)\delta-bk+1)^{2}}{b+1}\\
&+\frac{(-(2b^{3}+2b^{2}+b)\delta+2b^{3}k+b^{2}k+bk-3b^{2}-1)((2b+3)\delta-bk+1)}{b+1}\\
&+(b^{3}-b)\delta^{2}+\frac{(-2b^{4}k-2b^{3}k+bk+3b^{3}+3b^{2}+b-1)\delta}{b+1}\\
&+\frac{b^{4}k^{2}+b^{3}k^{2}+b^{2}k^{2}-3b^{3}k-2b^{2}k-3bk+2b^{2}+2}{b+1} \ \ \ \ \ (\mbox{since} \ n\geq(2b+3)\delta-bk+1)\\
=&\frac{(b^{4}+3b^{3}-4b)\delta^{2}+(2b^{3}k+6b^{2}k+4bk-b^{3}-2b^{2}-2b-4)\delta-b^{2}k-bk+1}{b+1}\\
>&0.
\end{align*}
Combining this with \eqref{eq:4.5}, \eqref{eq:4.6} and $s\geq\delta+1$, we obtain
$$
\varphi_{B_2}(\theta_*)=(s-\delta)f_1(\theta_*)>(s-\delta)f_1(n-b\delta+bk-2)=(s-\delta)f_2(s)>0.
$$
As $\theta_2<n-b\delta+bk-2<\rho(G_*)=\theta_*$ (see \eqref{eq:4.4}), we deduce that $\rho(G_2)=\theta_1<\theta_*=\rho(K_{\delta}\vee(K_{n-(b+1)\delta+bk-1}\cup(b\delta-bk+1)K_1))$. Together with \eqref{eq:4.1} and \eqref{eq:4.2},
we conclude
$$
\rho(G)\leq\rho(G_1)\leq\rho(G_2)<\rho(K_{\delta}\vee(K_{n-(b+1)\delta+bk-1}\cup(b\delta-bk+1)K_1)),
$$
which contradicts $\rho(G)\geq\rho(K_{\delta}\vee(K_{n-(b+1)\delta+bk-1}\cup(b\delta-bk+1)K_1))$.

\noindent{\bf Case 2.} $s=\delta$.

In terms of Lemma 2.6, we have
$$
\rho(G_1)\leq\rho(K_{\delta}\vee(K_{n-(b+1)\delta+bk-1}\cup(b\delta-bk+1)K_1)),
$$
where the equality holds if and only if $G_1=K_{\delta}\vee(K_{n-(b+1)\delta+bk-1}\cup(b\delta-bk+1)K_1)$. Combining this with \eqref{eq:4.1}, we
obtain
$$
\rho(G)\leq\rho(K_{\delta}\vee(K_{n-(b+1)\delta+bk-1}\cup(b\delta-bk+1)K_1)),
$$
with equality holding if and only if $G=K_{\delta}\vee(K_{n-(b+1)\delta+bk-1}\cup(b\delta-bk+1)K_1)$. Observe that $K_{\delta}\vee(K_{n-(b+1)\delta+bk-1}\cup(b\delta-bk+1)K_1)$ is not $k$-critical with respect to $[1,b]$-odd factor. Consequently, we can get a
contradiction.

\noindent{\bf Case 3.} $s\leq\delta-1$.

Let $G_3=K_s\vee(K_{n-s-(bs-bk+1)(\delta+1-s)}\cup(bs-bk+1)K_{\delta+1-s})$. Recall that the minimum degree of $G$ is $\delta$ and $G$ is a spanning
subgraph of $G_1=K_s\vee(K_{n_1}\cup K_{n_2}\cup\cdots\cup K_{n_{bs-bk+2}})$ for some positive odd integers $n_1\geq n_2\geq\cdots\geq n_{bs-bk+2}$
with $\sum\limits_{i=1}^{bs-bk+2}{n_i}=n-s$. Obviously, the minimum degree of $G_1$ is at least $\delta$, and so $n_{bs-bk+2}\geq\delta+1-s$.
According to Lemma 2.6, we conclude
\begin{align}\label{eq:4.7}
\rho(G_1)\leq\rho(G_3),
\end{align}
where the equality holds if and only if $(n_1,n_2,\ldots,n_{bs-bk+2})=(n-s-(bs-bk+1)(\delta+1-s),\delta+1-s,\ldots,\delta+1-s)$.

The quotient matrix of $A(G_3)$ corresponding to the partition $V(G_3)=V(K_s)\cup V(K_{n-s-(bs-bk+1)(\delta+1-s)})\cup V((bs-bk+1)K_{\delta+1-s})$
is given as
\begin{align*}
B_3=\left(
  \begin{array}{ccc}
    s-1 & n-s-(bs-bk+1)(\delta+1-s) & (bs-bk+1)(\delta+1-s)\\
    s & n-s-(bs-bk+1)(\delta+1-s)-1 & 0\\
    s & 0 & \delta-s\\
  \end{array}
\right).
\end{align*}
By a simple computation, the characteristic polynomial of $B_3$ is
\begin{align*}
\varphi_{B_3}(x)=&x^{3}+(-n-bs^{2}+b\delta s+bks+bs-bk\delta-bk+3)x^{2}\\
&+((b\delta-b)s^{2}+(bk+b+\delta-b\delta^{2}-bk\delta-n+1)s+\delta n-n+bk\delta^{2}-\delta^{2}-2\delta-bk+2)x\\
&+(-bs^{3}+b\delta s^{2}+bks^{2}+bs^{2}-s^{2}-bk\delta s-bks+\delta s+\delta)n-b^{2}s^{5}\\
&+(2b^{2}\delta+2b^{2}k+2b^{2}-b)s^{4}-(b^{2}\delta^{2}+4b^{2}k\delta+b^{2}k^{2}+2b^{2}\delta+4b^{2}k+b^{2}-3b\delta-bk-3b)s^{3}\\
&+(2b^{2}k\delta^{2}+2b^{2}k^{2}\delta+4b^{2}k\delta+2b^{2}k^{2}-3bk\delta+2b^{2}k-2b\delta^{2}-3b\delta-3bk-2b+\delta+1)s^{2}\\
&-(b^{2}k^{2}\delta^{2}+2b^{2}k^{2}\delta+b^{2}k^{2}-2bk\delta^{2}-3bk\delta+b\delta^{2}-2bk+b\delta+\delta^{2}+\delta)s\\
&+bk\delta^{2}+bk\delta-\delta^{2}-2\delta.
\end{align*}
Since the partition $V(G_3)=V(K_s)\cup V(K_{n-s-(bs-bk+1)(\delta+1-s)})\cup V((bs-bk+1)K_{\delta+1-s})$ is equitable, it follows from Lemma 2.7
that the largest root of $\varphi_{B_3}(x)=0$ equals $\rho(G_3)$.

By a direct calculation, the derivative function of $\varphi_{B_3}(x)$ is
\begin{align*}
\varphi_{B_3}'(x)=&3x^{2}+2(-n-bs^{2}+b\delta s+bks+bs-bk\delta-bk+3)x\\
&+(b\delta-b)s^{2}+(bk+b+\delta-b\delta^{2}-bk\delta-n+1)s+\delta n-n+bk\delta^{2}-\delta^{2}-2\delta-bk+2.
\end{align*}
Recall that $n\geq(bs-bk+2)(\delta+1-s)+s$. Note that
$$
-\frac{-n-bs^{2}+b\delta s+bks+bs-bk\delta-bk+3}{3}<n-b\delta+bk-2.
$$
Then $\varphi_{B_3}'(x)$ is increasing in the interval $[n-b\delta+bk-2,+\infty)$. For $x>n-b\delta+bk-2$, we deduce
\begin{align*}
\varphi_{B_3}'(x)>&\varphi_{B_3}'(n-b\delta+bk-2)\\
=&3(n-b\delta+bk-2)^{2}+2(-n-bs^{2}+b\delta s+bks+bs-bk\delta-bk+3)(n-b\delta+bk-2)\\
&+(b\delta-b)s^{2}+(bk+b+\delta-b\delta^{2}-bk\delta-n+1)s+\delta n-n+bk\delta^{2}-\delta^{2}-2\delta-bk+2\\
=&n^{2}+(-2bs^{2}+2b\delta s+2bks+2bs-s-2bk\delta-4b\delta+\delta+2bk-3)n\\
&+(2b^{2}\delta-2b^{2}k+b\delta+3b)s^{2}\\
&+(-2b^{2}\delta^{2}-b\delta^{2}-2b^{2}\delta-bk\delta-4b\delta+\delta+2b^{2}k^{2}+2b^{2}k-3bk-3b+1)s\\
&+2b^{2}k\delta^{2}+3b^{2}\delta^{2}+bk\delta^{2}-\delta^{2}-2b^{2}k^{2}\delta-4b^{2}k\delta+4bk\delta+6b\delta-2\delta+b^{2}k^{2}-3bk+2\\
\geq&(b\delta^{2}-bk)^{2}+(-2bs^{2}+2b\delta s+2bks+2bs-s-2bk\delta-4b\delta+\delta+2bk-3)(b\delta^{2}-bk)\\
&+(2b^{2}\delta-2b^{2}k+b\delta+3b)s^{2}\\
&+(-2b^{2}\delta^{2}-b\delta^{2}-2b^{2}\delta-bk\delta-4b\delta+\delta+2b^{2}k^{2}+2b^{2}k-3bk-3b+1)s\\
&+2b^{2}k\delta^{2}+3b^{2}\delta^{2}+bk\delta^{2}-\delta^{2}-2b^{2}k^{2}\delta-4b^{2}k\delta\\
&+4bk\delta+6b\delta-2\delta+b^{2}k^{2}-3bk+2 \ \ \ (\mbox{since} \ n\geq b\delta^{2}-bk)\\
=&\delta^{2}(b^{2}\delta^{2}+(2b^{2}s-2b^{2}k-4b^{2}+b)\delta-2b^{2}s^{2}+2b^{2}ks-2bs+2b^{2}k+3b^{2}+bk-3b-1)\\
&+\delta(2b^{2}s^{2}+bs^{2}-2b^{2}ks-2b^{2}s-bks-4bs+s+3bk+6b-2)\\
&+3bs^{2}-2bks-3bs+s+2\\
\geq&\delta^{2}(b^{2}(s+1)^{2}+(2b^{2}s-2b^{2}k-4b^{2}+b)(s+1)-2b^{2}s^{2}\\
&+2b^{2}ks-2bs+2b^{2}k+3b^{2}+bk-3b-1)\\
&+\delta(2b^{2}s^{2}+bs^{2}-2b^{2}ks-2b^{2}s-bks-4bs+s+3bk+6b-2)\\
&+3bs^{2}-2bks-3bs+s+2 \ \ \ (\mbox{since} \ \delta\geq s+1)\\
=&\delta^{2}(b^{2}s^{2}-bs+bk-2b-1)\\
&+\delta(2b^{2}s^{2}+bs^{2}-2b^{2}ks-2b^{2}s-bks-4bs+s+3bk+6b-2)\\
&+3bs^{2}-2bks-3bs+s+2\\
>&0 \ \ \ (\mbox{since} \ \delta\geq s+1\geq k+2\geq3),
\end{align*}
which implies that $\varphi_{B_3}(x)$ is increasing in the interval $[n-b\delta+bk-2,+\infty)$.

Recall that $G_*=K_{\delta}\vee(K_{n-(b+1)\delta+bk-1}\cup(b\delta-bk+1)K_1)$ and
\begin{align*}
\varphi_{B_*}(x)=&x^{3}+(-n+b\delta-bk+3)x^{2}+(-n-b\delta^{2}+bk\delta+b\delta-\delta-bk+2)x\\
&-b(b+1)\delta^{3}+bn\delta^{2}+2b^{2}k\delta^{2}+bk\delta^{2}-(3b+1)\delta^{2}-bkn\delta+n\delta-b^{2}k^{2}\delta+3bk\delta-2\delta.
\end{align*}
Notice that $\varphi_{B_*}(\theta_*)=0$. By plugging the value $\theta_*$ into $x$ of $\varphi_{B_3}(x)-\varphi_{B_*}(x)$, we get
\begin{align}\label{eq:4.8}
\varphi_{B_3}(\theta_*)=\varphi_{B_3}(\theta_*)-\varphi_{B_*}(\theta_*)=(\delta-s)g_1(\theta_*),
\end{align}
where $g_1(\theta_*)=(bs-bk-b)\theta_*^{2}+(n-b\delta s+bs+bk\delta+b\delta-bk-b-\delta-1)\theta_*+(bs^{2}-bks-bs-b\delta+bk+s)n+b^{2}s^{4}
-(b^{2}\delta+2b^{2}k+2b^{2}-b)s^{3}+(2b^{2}k\delta-2b\delta+b^{2}k^{2}+4b^{2}k+b^{2}-bk-3b)s^{2}-(b^{2}k^{2}\delta-b^{2}\delta-2bk\delta+2b^{2}k^{2}
+2b^{2}k-3bk-2b+\delta+1)s+b^{2}\delta^{2}+b\delta^{2}-2b^{2}k\delta+3b\delta+b^{2}k^{2}-2bk+2$. Note that $s\geq k+1$ by Claim 1. If $s=k+1$,
then $g_1(\theta_*)=(n-\delta-1)\theta_*+(-b\delta+bk+k+1)n+b^{2}\delta^{2}+b\delta^{2}-2b^{2}k\delta-2bk\delta+b\delta-k\delta-\delta
+b^{2}k^{2}+bk^{2}-bk-k+1$ and $n\geq(bs-bk+2)(\delta+1-s)+s=(b(k+1)-bk+2)(\delta+1-k-1)+k+1=b\delta+2\delta-bk-k+1$. Combining these with
\eqref{eq:4.4} and $k+1=s\leq\delta-1$, we obtain
\begin{align*}
g_1(\theta_*)>&(n-\delta-1)(n-b\delta+bk-2)+(-b\delta+bk+k+1)n+b^{2}\delta^{2}+b\delta^{2}\\
&-2b^{2}k\delta-2bk\delta+b\delta-k\delta-\delta+b^{2}k^{2}+bk^{2}-bk-k+1\\
=&n^{2}+(-2b\delta-\delta+2bk+k-2)n+b^{2}\delta^{2}+2b\delta^{2}-2b^{2}k\delta-3bk\delta\\
&+2b\delta-k\delta+\delta+b^{2}k^{2}+bk^{2}-2bk-k+3\\
\geq&(b\delta+2\delta-bk-k+1)^{2}+(-2b\delta-\delta+2bk+k-2)(b\delta+2\delta-bk-k+1)\\
&+b^{2}\delta^{2}+2b\delta^{2}-2b^{2}k\delta-3bk\delta+2b\delta-k\delta+\delta+b^{2}k^{2}+bk^{2}-2bk-k+3\\
=&(b+2)\delta^{2}-(bk+2k)\delta+2\\
\geq&(b+2)(k+2)\delta-(bk+2k)\delta+2\\
=&2(b+2)\delta+2\\
>&0.
\end{align*}
Combining this with \eqref{eq:4.8} and $s\leq\delta-1$, we infer
$$
\varphi_{B_3}(\theta_*)=(\delta-s)g_1(\theta_*)>0.
$$
Recall that $\varphi_{B_3}(x)$ is increasing in the interval $[n-b\delta+bk-2,+\infty)$. Then it follows from \eqref{eq:4.4} that
$$
\rho(G_3)<\theta_*=\rho(K_{\delta}\vee(K_{n-(b+1)\delta+bk-1}\cup(b\delta-bk+1)K_1)).
$$
Together with \eqref{eq:4.1} and \eqref{eq:4.7}, we obtain
$$
\rho(G)\leq\rho(G_1)\leq\rho(G_3)<\rho(K_{\delta}\vee(K_{n-(b+1)\delta+bk-1}\cup(b\delta-bk+1)K_1)),
$$
which is a contradiction to $\rho(G)\geq\rho(K_{\delta}\vee(K_{n-(b+1)\delta+bk-1}\cup(b\delta-bk+1)K_1))$. In what follows, we shall consider
$s\geq k+2$.

Recall that $g_1(\theta_*)=(bs-bk-b)\theta_*^{2}+(n-b\delta s+bs+bk\delta+b\delta-bk-b-\delta-1)\theta_*+(bs^{2}-bks-bs-b\delta+bk+s)n+b^{2}s^{4}
-(b^{2}\delta+2b^{2}k+2b^{2}-b)s^{3}+(2b^{2}k\delta-2b\delta+b^{2}k^{2}+4b^{2}k+b^{2}-bk-3b)s^{2}-(b^{2}k^{2}\delta-b^{2}\delta-2bk\delta+2b^{2}k^{2}
+2b^{2}k-3bk-2b+\delta+1)s+b^{2}\delta^{2}+b\delta^{2}-2b^{2}k\delta+3b\delta+b^{2}k^{2}-2bk+2$. Note that
$$
-\frac{n-b\delta s+bs+bk\delta+b\delta-bk-b-\delta-1}{2(bs-bk-b)}<n-b\delta+bk-2<\theta_*
$$
by \eqref{eq:4.4}, $n\geq(bs-bk+2)(\delta+1-s)+s$ and $k+2\leq s\leq\delta-1$. Then we have
\begin{align*}
g_1(\theta_*)>&g_1(n-b\delta+bk-2)\\
=&(bs-bk-b)(n-b\delta+bk-2)^{2}\\
&+(n-b\delta s+bs+bk\delta+b\delta-bk-b-\delta-1)(n-b\delta+bk-2)\\
&+(bs^{2}-bks-bs-b\delta+bk+s)n+b^{2}s^{4}-(b^{2}\delta+2b^{2}k+2b^{2}-b)s^{3}\\
&+(2b^{2}k\delta-2b\delta+b^{2}k^{2}+4b^{2}k+b^{2}-bk-3b)s^{2}\\
&-(b^{2}k^{2}\delta-b^{2}\delta-2bk\delta+2b^{2}k^{2}+2b^{2}k-3bk-2b+\delta+1)s\\
&+b^{2}\delta^{2}+b\delta^{2}-2b^{2}k\delta+3b\delta+b^{2}k^{2}-2bk+2\\
=&(bs-bk-b+1)n^{2}+(bs^{2}-2b^{2}\delta s+2b^{2}ks-b\delta s-bks-4bs+s\\
&+2b^{2}k\delta+2b^{2}\delta+bk\delta-b\delta-2b^{2}k^{2}-2b^{2}k+5bk+3b-\delta-3)n\\
&+b^{2}s^{4}-(b^{2}\delta+2b^{2}k+2b^{2}-b)s^{3}\\
&+(2b^{2}k\delta-2b\delta+b^{2}k^{2}+4b^{2}k+b^{2}-bk-3b)s^{2}\\
&+(b^{3}\delta^{2}+b^{2}\delta^{2}-2b^{3}k\delta-b^{2}k^{2}\delta-b^{2}k\delta+4b^{2}\delta+2bk\delta+2b\delta-\delta\\
&+b^{3}k^{2}-2b^{2}k^{2}-5b^{2}k+3bk+4b-1)s\\
&-b^{3}k\delta^{2}-b^{3}\delta^{2}-b^{2}k\delta^{2}+2b\delta^{2}+2b^{3}k^{2}\delta+2b^{3}k\delta+b^{2}k^{2}\delta-4b^{2}k\delta\\
&-3b^{2}\delta-3bk\delta+2b\delta+2\delta-b^{3}k^{3}-b^{3}k^{2}+4b^{2}k^{2}+3b^{2}k-5bk-2b+4.
\end{align*}
Notice that $-(bs^{2}-2b^{2}\delta s+2b^{2}ks-b\delta s-bks-4bs+s+2b^{2}k\delta+2b^{2}\delta+bk\delta-b\delta-2b^{2}k^{2}-2b^{2}k+5bk+3b-\delta-3)/(2(bs-bk-b+1))
=-\frac{s-2b\delta+2bk-\delta-3}{2}\leq\frac{2b\delta+\delta-2bk}{2}<b\delta^{2}-bk\leq n$, and so
\begin{align*}
g_1(\theta_*)>&g_1(n-b\delta+bk-2)\\
\geq&(bs-bk-b+1)(b\delta^{2}-bk)^{2}+(bs^{2}-2b^{2}\delta s+2b^{2}ks-b\delta s-bks-4bs+s\\
&+2b^{2}k\delta+2b^{2}\delta+bk\delta-b\delta-2b^{2}k^{2}-2b^{2}k+5bk+3b-\delta-3)(b\delta^{2}-bk)\\
&+b^{2}s^{4}-(b^{2}\delta+2b^{2}k+2b^{2}-b)s^{3}\\
&+(2b^{2}k\delta-2b\delta+b^{2}k^{2}+4b^{2}k+b^{2}-bk-3b)s^{2}\\
&+(b^{3}\delta^{2}+b^{2}\delta^{2}-2b^{3}k\delta-b^{2}k^{2}\delta-b^{2}k\delta+4b^{2}\delta+2bk\delta+2b\delta-\delta\\
&+b^{3}k^{2}-2b^{2}k^{2}-5b^{2}k+3bk+4b-1)s\\
&-b^{3}k\delta^{2}-b^{3}\delta^{2}-b^{2}k\delta^{2}+2b\delta^{2}+2b^{3}k^{2}\delta+2b^{3}k\delta+b^{2}k^{2}\delta-4b^{2}k\delta\\
&-3b^{2}\delta-3bk\delta+2b\delta+2\delta-b^{3}k^{3}-b^{3}k^{2}+4b^{2}k^{2}+3b^{2}k-5bk-2b+4\\
=&\delta(b\delta((b^{2}s-b^{2}k-b^{2}+b)\delta^{2}-(2b^{2}s+bs-2b^{2}k-2b^{2}-bk+b+1)\delta\\
&+bs^{2}+b^{2}s-bks-3bs+s-b^{2}k+2bk-b^{2}+3b-1)-b^{2}s^{3}+2b^{2}ks^{2}-2bs^{2}\\
&-b^{2}k^{2}s+4b^{2}s+2bks+2bs-s-3b^{2}k-3b^{2}-2bk+2b+2)+b^{2}s^{4}\\
&-(2b^{2}k+2b^{2}-b)s^{3}+(b^{2}k^{2}+3b^{2}k+b^{2}-bk-3b)s^{2}\\
&-(b^{2}k^{2}+b^{2}k-2bk-4b+1)s-2bk-2b+4\\
\geq&\delta(b\delta((b^{2}s-b^{2}k-b^{2}+b)(s+1)^{2}-(2b^{2}s+bs-2b^{2}k-2b^{2}-bk+b+1)(s+1)\\
&+bs^{2}+b^{2}s-bks-3bs+s-b^{2}k+2bk-b^{2}+3b-1)-b^{2}s^{3}+2b^{2}ks^{2}-2bs^{2}\\
&-b^{2}k^{2}s+4b^{2}s+2bks+2bs-s-3b^{2}k-3b^{2}-2bk+2b+2)+b^{2}s^{4}\\
&-(2b^{2}k+2b^{2}-b)s^{3}+(b^{2}k^{2}+3b^{2}k+b^{2}-bk-3b)s^{2}\\
&-(b^{2}k^{2}+b^{2}k-2bk-4b+1)s-2bk-2b+4 \ \ \ \ \ (\mbox{since} \ \delta\geq s+1\geq k+3)\\
=&\delta(b\delta(b^{2}s^{3}-b^{2}ks^{2}-b^{2}s^{2}+bs^{2}-3bs+3bk+3b-2)-b^{2}s^{3}+2b^{2}ks^{2}-2bs^{2}\\
&-b^{2}k^{2}s+4b^{2}s+2bks+2bs-s-3b^{2}k-3b^{2}-2bk+2b+2)+b^{2}s^{4}\\
&-(2b^{2}k+2b^{2}-b)s^{3}+(b^{2}k^{2}+3b^{2}k+b^{2}-bk-3b)s^{2}\\
&-(b^{2}k^{2}+b^{2}k-2bk-4b+1)s-2bk-2b+4\\
\geq&\delta(b(s+1)(b^{2}s^{3}-b^{2}ks^{2}-b^{2}s^{2}+bs^{2}-3bs+3bk+3b-2)-b^{2}s^{3}+2b^{2}ks^{2}\\
&-2bs^{2}-b^{2}k^{2}s+4b^{2}s+2bks+2bs-s-3b^{2}k-3b^{2}-2bk+2b+2)+b^{2}s^{4}\\
&-(2b^{2}k+2b^{2}-b)s^{3}+(b^{2}k^{2}+3b^{2}k+b^{2}-bk-3b)s^{2}\\
&-(b^{2}k^{2}+b^{2}k-2bk-4b+1)s-2bk-2b+4 \ \ \ \ \ (\mbox{since} \ \delta\geq s+1\geq k+3)\\
=&\delta(b^{3}s^{4}-b^{3}ks^{3}-b^{3}ks^{2}-b^{3}s^{2}+2b^{2}ks^{2}-2b^{2}s^{2}-2bs^{2}-b^{2}k^{2}s+3b^{2}ks+4b^{2}s\\
&+2bks-s-2bk+2)+b^{2}s^{4}-(2b^{2}k+2b^{2}-b)s^{3}+(b^{2}k^{2}+3b^{2}k+b^{2}-bk-3b)s^{2}\\
&-(b^{2}k^{2}+b^{2}k-2bk-4b+1)s-2bk-2b+4\\
\geq&(k+3)(b^{3}s^{4}-b^{3}ks^{3}-b^{3}ks^{2}-b^{3}s^{2}+2b^{2}ks^{2}-2b^{2}s^{2}-2bs^{2}-b^{2}k^{2}s+3b^{2}ks+4b^{2}s\\
&+2bks-s-2bk+2)+b^{2}s^{4}-(2b^{2}k+2b^{2}-b)s^{3}+(b^{2}k^{2}+3b^{2}k+b^{2}-bk-3b)s^{2}\\
&-(b^{2}k^{2}+b^{2}k-2bk-4b+1)s-2bk-2b+4 \ \ \ \ \ (\mbox{since} \ \delta\geq s+1\geq k+3)\\
=&(b^{3}k+3b^{3}+b^{2})s^{4}-(b^{3}k^{2}+3b^{3}k+2b^{2}k+2b^{2}-b)s^{3}\\
&-(b^{3}k^{2}+4b^{3}k-3b^{2}k^{2}-7b^{2}k+3bk+3b^{3}+5b^{2}+9b)s^{2}\\
&-(b^{2}k^{3}+b^{2}k^{2}-12b^{2}k-2bk^{2}-8bk+k-12b^{2}-4b+4)s\\
&-2bk^{2}-8bk-2b+2k+10\\
>&0 \ \ \ \ \ (\mbox{since} \ s\geq k+2).
\end{align*}
Combining this with \eqref{eq:4.8} and $s\leq\delta-1$, we conclude
\begin{align}\label{eq:4.9}
\varphi_{B_3}(\theta_*)=(\delta-s)g_1(\theta_*)>0.
\end{align}
Recall that $\varphi_{B_3}(x)$ is increasing in the interval $[n-b\delta+bk-2,+\infty)$. Together with \eqref{eq:4.4} and \eqref{eq:4.9}, we get
$$
\rho(G_3)<\theta_*=\rho(K_{\delta}\vee(K_{n-(b+1)\delta+bk-1}\cup(b\delta-bk+1)K_1)).
$$
Combining this with \eqref{eq:4.1} and \eqref{eq:4.7}, we conclude
$$
\rho(G)\leq\rho(G_1)\leq\rho(G_3)<\rho(K_{\delta}\vee(K_{n-(b+1)\delta+bk-1}\cup(b\delta-bk+1)K_1)),
$$
which contradicts $\rho(G)\geq\rho(K_{\delta}\vee(K_{n-(b+1)\delta+bk-1}\cup(b\delta-bk+1)K_1))$. This completes the proof of Theorem 1.2. \hfill $\Box$

\medskip

\section*{Data availability statement}

My manuscript has no associated data.

\section*{Declaration of competing interest}

The authors declare that they have no conflicts of interest to this work.

\section*{Acknowledgments}

This work was supported by the Natural Science Foundation of Jiangsu Province (Grant No. BK20241949). Project ZR2023MA078 supported by Shandong
Provincial Natural Science Foundation.

\end{document}